\documentclass[11pt]{amsart}

\usepackage{latexsym}
\usepackage{amssymb}
\usepackage{amsfonts}
\usepackage{fullpage}

\newtheorem{thm}{Theorem}[section]
\newtheorem{lem}[thm]{Lemma}
\newtheorem{cor}[thm]{Corollary}
\newtheorem{pro}[thm]{Proposition}

\theoremstyle{definition}
\newtheorem{defn}[thm]{Definition}
\newtheorem{exmp}[thm]{Example}
\newtheorem{conj}[thm]{Conjecture}
\newtheorem{rem}[thm]{Remark}
\newtheorem{que}[thm]{Question}

\def\H{\text{\bf H}}
\def\P{{\mathbb P}}

\def\X{{\mathbb X}}

\def\T{{\mathcal T}}
\def\H{{\bf H}}

\def\ni{\noindent}

\def\ra{\rightarrow}

\def\ds{\displaystyle}

\def\FF{{\mathcal F}}

\def\GG{{\mathcal G}}

\def\proj{\mathbb P}

\begin{document}

\title[Non-Level O-sequences]{Non-Level O-sequences of Codimension $3$ and  Degree of The Socle Elements}
 
\date{\today}

\author[Y.S.\, Shin]{Yong Su Shin}
\address{Department of Mathematics,  Sungshin Women's University,
Seoul,  136-742, Korea}
\email{ysshin@sungshin.ac.kr}
\thanks{This work was supported by Korea Research Foundation Grant        (KRF-2003-015-C00004).}
\begin{abstract} 
It is unknown if an Artinian level O-sequence of codimension $3$ and type $r\ (\ge 2)$ is unimodal, while it is known that any Gorenstein O-sequence of codimension $3$ is unimodal. We show that some Artinian non-unimodal O-sequence of codimension $3$ cannot be level. 
We also  find another non-level case: if some Artinian algebra $A$ of codimension $3$ has the Hilbert function 
$$
\begin{matrix}\H & : & h_0 & h_1 & \cdots & h_{d-1} & \underbrace{h_d \ \cdots \ h_d}_{s\text{-times}} & h_{d+s}, \end{matrix}
$$
such that $h_d<h_{d+s}$ and $s\ge 2$, then $A$ has a socle element in degree $d+s-2$, that is, $A$ is not level. 
\end{abstract}
\maketitle 

\section{Introduction}

Let ${\mathbb X} = \{ P_1, \ldots , P_s \}$ be a set of $s$ distinct points in the projective space $\proj ^n(k)$ (where $k = \overline{k}$ is an algebraically closed field).  Then $P_i \leftrightarrow \wp _i = (L_{i1}, \ldots , L_{in}) \subset R = k[x_0, x_1, \ldots , x_n]$ where the $L_{ij}, \ j = 1, \ldots , n$ are $n$ linearly independent linear forms and $\wp_i$ is the (homogeneous) prime ideal of $R$ generated by all the forms which vanish at $P_i$.  The ideal
$$
I = I_{\mathbb X}:= \wp _1\cap \cdots \cap \wp_s
$$
is the ideal generated by all the forms which vanish at all the points of $\mathbb X$.

Since $R = \oplus_{i=0}^\infty R_i$ ($R_i$: the vector space of dimension $\binom{i+n}{n}$ generated by all the monomials in $R$ having degree $i$) and $I = \oplus_{i=0}^\infty I_i$, we get that 
$$
A = R/I = \oplus_{i=0}^\infty (R_i/I_i) = \oplus _{i=0}^\infty A_i
$$
is a graded ring.  The numerical function
$$
\H_{\mathbb X}(t) = \H_A(t):= \dim _k A_t = \dim _kR_t - \dim _k I_t
$$
is called the {\em Hilbert function} of the set $\mathbb X$ (or of the ring $A$). 

Given an O-sequence $\H=(h_0,h_1,\dots)$, we define the {\em first difference} of $\H$ as 
$$
\Delta \H=(h_0,h_1-h_0,h_2-h_1,h_3-h_2,\dots). 
$$

Let $h$ and $i$ be positive integers. Then $h$ can be written uniquely in the form
$$
h=\binom{m_i}{i}+\binom{m_{i-1}}{i-1}+\cdots+\binom{m_j}{j}
$$
where $m_i>m_{i-1}>\cdots >m_j\ge j\ge 1$. This expansion for $h$ is called the $i$-{\em binomial expansion} of $h$. Also, define
$$
h^{\langle i\rangle}=\binom{m_i+1}{i+1}+\binom{m_{i-1}+1}{(i-1)+1}+\cdots+\binom{m_j+1}{j+1},
$$
and $0^{\langle i\rangle}=0$. 

It is worth noting that $R$ is a {\em standard} graded algebra 
since $R=k[R_1]$, that is, $R$ is generated (as a $k$-algebra) by its piece of degree $1$. If $I$ is a homogeneous ideal of $R$, then $R/I$ is again a standard graded $k$-algebra. Furthermore, if $I$ has a height $n+1$ in $R$, then $A=R/I$ is an {\em Artinian} $k$-algebra, and hence $\dim_k A<\infty$. Thus we can write $A=k\oplus A_1\oplus \cdots \oplus A_s$ where $A_s\ne 0$. We call $s$ the {\em socle degree} of $A$.

We associate to the graded Artinian algebra $A$ a vector  of non-negative integers which is an $(s+1)$-tuple, called the {\em $h$-vector} of $A$ and denoted by $h(A)$. It is defined as follows.
$$
h(A):=(1,\dim_k A_1,\dots,\dim_k A_s)=(h_0,h_1,\dots,h_s)
\ \ \ {\text{with}\ \ \ } h_s\ne 0. 
$$
Moreover, an $h$-vector $(h_0,h_1,\dots,h_s)$ is called {\em unimodal} if $h_0\le \cdots \le h_t=\cdots=h_\ell \ge \cdots \ge h_s$. 

Let $\FF_\X$ be the graded minimal resolution of $R/I_\X$ (or $\X$), i.e.,
$$
\begin{matrix}
\FF_\X: & 0 & \ra & \FF_n & \ra & \FF_{n-1} & \ra & \cdots
     & \ra & \FF_1 & \ra & R & \ra & R/I_\X & \ra & 0.
\end{matrix}
$$
We can write 
$$
\FF_i
 =  \bigoplus^{\gamma_{i}}_{j=1} R^{\beta_{ij}}(-\alpha_{ij})
$$ 
where $\alpha_{i1}< \alpha_{i2}< \cdots < \alpha_{i\gamma_{i}}$. The numbers $\alpha_{ij}$ are called the {\em shifts} associated to $R/I_\X$, and the numbers $\beta_{ij}$ are called the {\em graded Betti numbers} of $R/I_\X$ (or $\X$). 

Now, we recall that if the last free module of the minimal free resolution of a graded ring $A$ with Hilbert function $\H$ is of the form $\FF_n=R^\beta(-s)$ for some $s>0$, then   Hilbert function $\H$ and a graded ring $A$ are called {\em level}. In particular, if $\FF_n=R^\beta(-s)$ in $\FF_\X$, then we call $\X$ a {\em  level set of points} in $\P^n$.  For a special case, if $\beta=1$, then we call a graded Artinian algebra $A$ {\em Gorenstein}. 
 
In~\cite{St}, Stanley proved that any graded Artinian Gorenstein algebra of codimension $3$ is unimodal. In fact, he proved a stronger result than unimodality using the structure theorem of Buchsbaum and Eisenbud for the Gorenstein algebra of codimension $3$ in \cite{BE}. Since then, the graded Artinian Gorenstein algebras of codimension $3$ have been much studied (see \cite{D}, \cite{GHMS}, \cite{GHS:2},  \cite{Ha:1}, \cite{Ha:2}, \cite{Sh:1}). In~\cite{BI}, Bernstein and Iarrobino showed how to construct non-unimodal graded Artinian Gorenstein  algebras of codimension higher than or equal to $5$. Moreover, in \cite{BL}, Boij and Laksov showed another method on how to construct the same graded Artinian Gorenstein  algebras. Unfortunately, it has been unknown if there exists a graded non-unimodal Gorenstein algebra of codimension $4$. For unimodal Artinian Gorenstein algebras of codimension $4$, it has been shown in \cite{Sh:1} how to construct  some of them using the link-sum method. We have also shown in \cite{GHS:2} how to obtain some of unimodal Artinian Gorenstein algebras of any codimension $n\ (\ge 3)$. 

For graded Artinian level algebras, it has  been recently studied (see, \cite{BI},  \cite{BL}, \cite{GHMS}, \cite{GHS:4}). Since every graded Artinian Gorenstein algebra of codimension $3$ is unimodal, the following question in~\cite{GHMS} is quite interesting.

\begin{que} \label{Q:051}
Is any level O-sequence of codimension $3$ unimodal (Question 4.4, \cite{GHMS})?
\end{que} 

In~\cite{GHMS}, we proved the following result.  Let 
\begin{equation} \label{EQ:011} 
\begin{array}{llllllllllllllllllllll}
\H & : & h_0 & h_1 & \cdots & h_{d-1} & h_d & h_d & \cdots 
\end{array}
\end{equation} 
with $h_{d-1}>h_d$. If  $h_d\le d+1$ with any codimension $h_1$, then $\H$ is not {\em level} (see Proposition~\ref{P:006}). 

The goal of this paper is to find an answer to Question~\ref{Q:051} and we give an answer to this question under a certain condition.  In fact, it suffices to find an answer to the following Question~\ref{Q:012} (see Corollary~\ref{C:510}). 

\begin{que} \label{Q:012} 
Let $\H$ be an O-sequence as in equation~\eqref{EQ:011} with codimension $3$. Is $\H$ NOT level? 
\end{que}

As we mentioned above, it is shown that any Hilbert function $\H$ in equation~\eqref{EQ:011} is not level when $h_d\le d+1$. In Section~\ref{Sec:005}, we prove that any Hilbert function $\H$ with codimension $3$ in equation \eqref{EQ:011} is not level when $h_d\le 2d+2$ (see Theorem~\ref{T:014}). This provides an answer to Question~\ref{Q:051} when $h_d\le 2d+2$ (see Corollary~\ref{C:510}). Finally, in Section~\ref{Sec:006}, we find the degree of the socle elements of a graded Artinian algebra of codimension $3$ with Hilbert function
$$
\begin{matrix}
\H & : & h_0 & h_1 & \cdots & h_{d-1} & \underbrace{h_d \ \cdots \ h_d}_{s\text{-times}} & h_{d+s},
\end{matrix}
$$
where $h_d<h_{d+s}$ and $s\ge 2$. We prove that some graded algebra with Hilbert function $\H$ is not level and has a socle element in degree $d+s-2$ (see Theorem~\ref{T:062}). 

\section{Some Non-Level O-sequences of Codimension $3$} \label{Sec:005} 

In~\cite{GHMS}, we got an answer to {Question~\ref{Q:051}} with the condition $h_d\le d+1$ as follows. 

\begin{pro}[Proposition 2.21, {\cite{GHMS}}] \label{P:006} Let $h = (1, n, h_2, \ldots , h_s)$ be the $h$-vector of an Artinian algebra with socle degree $s$.  Then $h$ is {\bf\em not} a level sequence if $h_d=h_{d+1}\le d+1$ and $h_{d-1}>h_d$.
\end{pro} 

We shall expand the above proposition to a case $h_d\le 2d+2$ with codimension $3$. 

Let $R=k[x_1,\dots,x_n]$ and $A=R/I$ where $I$ is a homogeneous ideal of $R$ having height $n$. Then $A$ has the minimal free resolution $\FF$, as an $R$-module, of the form:
$$
\begin{array}{lllllllllllllll}
0 & \ra & \FF_{n-1} & \ra & \cdots & \ra & \FF_1 & \ra & \FF_0 & \ra & R & \ra
& A & \ra & 0
\end{array}
$$ 
where $\FF_j=\bigoplus_{t=1}^{\gamma_j} R^{\beta_{j,j+1+t}}(-(j+1+t))$  are each free graded $R$-modules. In \cite{EK}, Eliahou and Kervaire studied minimal free resolutions of certain monomial ideals. We recall some of their notations and   results  here.

\begin{defn} \label{D:051}
Let $T\in R=k[x_1,\dots,x_n]$ be a {\em term} of $R$. Then
$$
m(T):=\max\{\,i \mid x_i \text{ divides } T\, \}.
$$
In other words, $m(T)$ is the largest index of an indeterminate that divides $T$.
\end{defn} 

\begin{thm} [Eliahou--Kervaire, {\cite{EK}}] \label{T:052}
Let $I$ be a stable monomial ideal of $R$ (e.g., a lex segment ideal). Denote by $\GG(I)_d$ the elements of that set which have degree $d$. Then
$$
\beta_{q,i} = \sum_{T\in \GG(I)_{i-q}} \binom{m(T)-1}{q}.
$$
\end{thm}

This beautiful theorem gives all the graded Betti numbers of the lex segment ideal just from an intimate knowledge of the generators of that ideal. Since the minimal free resolution of the ideal of a $k$-configuration in $\P^n$ is extremal (\cite{GHS:2}, \cite{GS:1}), we may apply this result to those ideals. It is an immediate consequence of the Eliahou--Kervaire theorem that if $I$ is either a lex-segment ideal or the ideal of a $k$-configuration in $\P^n$ which has {\em no} generators of degree $d$, then $\beta_{q,i}=0$ whenever $i-q=d$. 

By the result of \cite{Pe}, the only way we can cancel graded Betti numbers is if there are the same graded Betti numbers in the adjacent free modules of the extremal minimal free resolution. Note that it is quite obvious for a case of $n=3$.  

The following lemma is a simple consequence of a lex segment ideal, so we shall omit the proof here. 

\begin{lem} \label{L:013}
Let $I$ be the  lex-segment ideal in $R=k[x_1,x_2,x_3]$ with Hilbert function $\H=(h_0,h_1,\dots,h_s)$ where  $h_d=d+i$ and $1\le i\le \frac{d^2+d}{2}$.  Then the last monomial of $I_d$ is
$$
\begin{array}{clrlllllllllll}
x_1x_2^{i-1}x_3^{d-i}, & \text{ for } & 1 & \le & i & \le & d, \\
x_1^2x_2^{i-(d+1)}x_3^{(2d-1)-i}, & \text{ for } & d+1 & \le & i & \le & 2d-1,\\
\vdots\\
x_1^{d-1}x_2^{i-\frac{d^2+d-4}{2}}x_3^{\frac{d^2+d-2}{2}-i}, & \text{ for } & \frac{d^2+d-4}{2} & \le & i & \le & \frac{d^2+d-2}{2},\\
x_1^{d}, & \text{ for } &   &   & i & = & \frac{d^2+d}{2}.
\end{array}
$$
\end{lem}

We need the following proposition to prove the main Theorem~\ref{T:014}.

\begin{pro} \label{P:014}
Let $R=k[x_1,x_2,x_3]$ and let $\H=(h_0,h_1,\dots,h_s)$ be the $h$-vector of an Artinian algebra with socle degree $s$ and 
$$
h_d=h_{d+1}=d+i, \quad h_{d-1}>h_d, \quad \text{and} \quad j:=h_{d-1}-h_d
$$
for $i=1,2,\dots,\frac{d^2+d}{2}$. Then, for every $1\le k\le d$ and $1\le \ell \le d$, 
$$
\begin{array}{llllllllll}
\ds \beta_{1,d+2}
= \begin{cases}
2k-1, & \text{for } \quad 
(k-1)d-\frac{k(k-3)}{2} \le i \le (k-1)d-\frac{k(k-3)}{2}+(k-1),
\\
2k, & \text{for } \quad (k-1)d-\frac{k(k-3)}{2}+k\le i \le kd-\frac{(k-1)k}{2}.
\end{cases}\\[5ex]
{\ds \beta_{2,d+2}} 
= j+\ell, \hskip 7.9 mm  
\text{ for } \quad (\ell-1) d-\frac{(\ell-2)(\ell-1)}{2}<i\le  \ell d-\frac{(\ell-1)\ell}{2}.
\end{array}
$$
\end{pro}

\begin{proof} 
Since we assume $h_d=d+i$, the monomials not in $I_d$ are the last $d+i$ monomials of $R_d$. 
By Lemma~\ref{L:013}, the last monomial of $R_1I_d$ is
$$
\begin{array}{clllllllllllll}
x_1x_2^{i-1}x_3^{d-i+1}, & \text{ for } & i=1,\dots, d, \\
x_1^2x_2^{i-(d+1)}x_3^{2d-i}, & \text{ for }  & i=d+1,\dots,2d-1,\\
\vdots\\
x_1^{d-1}x_2^{i-\frac{d^2+d-4}{2}}x_3^{\frac{d^2+d}{2}-i}, & \text{ for } &  i=\frac{d^2+d-4}{2},\ \frac{d^2+d-2}{2}, \\
x_1^{d}x_3, & \text{ for } &    i = \frac{d^2+d}{2}.
\end{array}
$$
In what follows, the first monomial of $I_{d+1}-R_1I_d$ is 
\begin{equation}\label{EQ:001}
\begin{array}{cllllllllllll}
x_2^{d+1}, & \text{ for } & i=1, \\
x_1x_2^{i-2}x_3^{(d+2)-i}, & \text{ for } & i=2,\dots,d, \\
\vdots\\
x_1^{d-1}x_2x_3, & \text{ for } & i=\frac{d^2+d-2}{2}, \\
x_1^{d-1}x_2^2, & \text{ for } &   i=\frac{d^2+d}{2}. 
\end{array}
\end{equation}
Note that
\begin{equation}\label{EQ:002}
\begin{array}{lllllllllllll}
(d+i)^{\langle d\rangle}
& = & (d+i)+k, & \text{ for }\ i=(k-1)d-\frac{k(k-3)}{2},\dots, 
kd-\frac{k(k-1)}{2},\\
&   &          & \text{ and }k=1,\dots,d.
\end{array}
\end{equation}

We now calculate the Betti number
$$
\ds \beta_{1,d+2}
= \ds\sum_{T\in \GG(I)_{d+1}} \binom{m(T)-1}{1}.
$$
Based on equation~\eqref{EQ:001}, we shall find this Betti number of two cases for $i$ as follows. 
\begin{enumerate}

\item[]{\em Case 1-1.}  $i=(k-1)d-\frac{k(k-3)}{2}$ and $k=1,2,\dots,d$. 

Then, by equation~\eqref{EQ:002}, $I_{d+1}$ has $k$-generators, which are
$$
x_1^{k-1} x_2^{(d+2)-k},x_1^{k-1} x_2^{(d+1)-k}x_3,\dots,x_1^{k-1} x_2^{(d+3)-2k} x_3^{k-1}. 
$$
By the similar argument, for $i=(k-1)d-\frac{k(k-3)}{2}+1,\dots, (k-1)d-\frac{k(k-3)}{2}+(k-1)$, $I_{d+1}$ has $k$-generators including the element $x_1^{k-1} x_2^{(d+2)-k}$. Hence we have that 
$$
\ds \beta_{1,d+2}
= \ds\sum_{T\in \GG(I)_{d+1}} \binom{m(T)-1}{1}=2\times (k-1)+1=2k-1. 
$$

\item[]{\em Case 1-2.} $i=(k-1)d-\frac{k(k-3)}{2}+k=(k-1)d-\frac{k(k-5)}{2},\dots,
kd-\frac{k(k-1)}{2}$ and 
$k=1,2,\dots,d$. 

Then, by equation~\eqref{EQ:002}, $I_{d+1}$ has $k$-generators, which are
$$
x_1^{k} x_2^{i-\left((k-1)d-\frac{k^2-3k-2}{2}\right)}
x_3^{kd-\frac{k^2-k-4}{2}-i},\dots, 
x_1^{k} x_2^{i-\left((k-1)d-\frac{k(k-5)}{2}\right)}
x_3^{\left(kd-\frac{k(k-3)}{2}+1\right)-i}.
$$
Hence we have that 
$$
\ds \beta_{1,d+2}
= \ds\sum_{T\in \GG(I)_{d+1}} \binom{m(T)-1}{1}=2\times k=2k. 
$$
\end{enumerate}
 
Now we move on the Betti number:
$$
\ds \beta_{2,d+2}=\sum_{T\in \GG(I)_d} \binom{m(T)-1}{2}.
$$
Recall $h_d=d+i$ and $j:=h_{d-1}-h_d$. The calculation in this case is much more complicated, and there are four cases based on $i$ and $j$.  

\begin{enumerate}

\item[]{\em Case 2-1.} $(\ell-1) d-\frac{(\ell-2)(\ell-1)}{2}<i< \ell d-\frac{(\ell-1)\ell}{2}$ and $\ell=1,2,\dots,d$. 

Then the last monomial of $I_d$ is
$$
x_1^{\ell} x_2^{i-(\ell-1) d+\frac{\ell(\ell-3)}{2}} x_3^{\ell d -\frac{(\ell-1)\ell}{2}-i}.  
$$ 

\begin{enumerate}

\item  $(k-1)d-\frac{(k-1)k}{2}<i+j<kd-\frac{k(k+1)}{2}$ and $k=\ell,\ell+1,\dots,d$.

Then the first monomial of $I_d-R_1I_{d-1}$ is
$$
x_1^{k} x_2^{(i+j)-\left((k-1)d-\frac{(k-2)(k+1)}{2}\right)}x_3^{\left(kd-\frac{(k-1)(k+2)}{2}\right)-(i+j)},
$$
and hence we have $(j+k)$-generators in $I_d$ as follows:
$$
\begin{array}{llllllllllllllll}
x_1^{k} x_2^{(i+j)-\left((k-1)d-\frac{(k-2)(k+1)}{2}\right)}x_3^{\left(kd-\frac{(k-1)(k+2)}{2}\right)-(i+j)}, \dots,x_1^{k}x_3^{d-k},\\
x_1^{(k-1)} x_2^{d-(k-1)}, x_1^{(k-1)} x_2^{(d-1)-(k-1)}x_3,\dots, x_1^{(k-1)} x_3^{d-(k-1)}, \\
\hskip 2.5 true cm \vdots \\

x_1^{\ell+1} x_2^{(d-1)-\ell}, x_1^{\ell+1} x_2^{(d-2)-\ell}x_3,
\dots, x_1^{\ell+1} x_3^{(d-1)-\ell}\\
x_1^{\ell} x_2^{d-\ell},  \dots,
x_1^{\ell} x_2^{i-(\ell-1) d+\frac{\ell(\ell-3)}{2}} x_3^{\ell d -\frac{(\ell-1)\ell}{2}-i}\\
\end{array}
$$
and thus
$$
\ds \beta_{2,d+2}=\sum_{T\in \GG(I)_d} \binom{m(T)-1}{2}=j+\ell,
$$

\item $i+j=(k-1)d-\frac{(k-1)k}{2}$ and $k=\ell+1,\dots,d$.

Then the first monomial of $I_d-R_1I_{d-1}$ is
$$
x_1^{k-1} x_2^{d-(k-1)},
$$
and hence we have $(j+k)$-generators in $I_d$ as follows:
$$
\begin{array}{llllllllllllllll}
x_1^{k-1} x_2^{d-(k-1)}, x_1^{k-1} x_2^{(d-1)-(k-1)}x_3,\dots, x_1^{k-1} x_3^{d-(k-1)}, \\
\hskip 2.5 true cm \vdots \\
x_1^{\ell+1} x_2^{(d-1)-\ell}, x_1^{\ell+1} x_2^{(d-2)-\ell}x_3,
\dots, x_1^{\ell+1} x_3^{(d-1)-\ell}\\
x_1^{\ell} x_2^{d-\ell},  \dots,
x_1^{\ell} x_2^{i-(\ell-1) d+\frac{\ell(\ell-3)}{2}} x_3^{\ell d -\frac{(\ell-1)\ell}{2}-i}\\
\end{array}
$$
and thus
$$
\ds \beta_{2,d+2}=\sum_{T\in \GG(I)_d} \binom{m(T)-1}{2}=j+\ell.
$$

\end{enumerate}

\item[]{\em Case 2-2.} $i=\ell d-\frac{(\ell-1)\ell}{2}$ and $\ell=1,2,\dots,d$. 

Then the last monomial of $I_d$ is
$$
x_1^\ell x_2^{d-\ell}.  
$$

\begin{enumerate}

\item $(k-1)d-\frac{(k-1)k}{2}<i+j<kd-\frac{k(k+1)}{2}$ and $k=\ell+1,\dots,d$.

Then the first monomial of $I_d-R_1I_{d-1}$ is
$$
x_1^{k} x_2^{(i+j)-\left((k-1)d-\frac{(k-2)(k+1)}{2}\right)}x_3^{\left(kd-\frac{(k-1)(k+2)}{2}\right)-(i+j)},
$$
and hence we have $(j+k)$-generators in $I_d$ as follows:
$$
\begin{array}{llllllllllllllll}
x_1^{k} x_2^{(i+j)-\left((k-1)d-\frac{(k-2)(k+1)}{2}\right)}x_3^{\left(kd-\frac{(k-1)(k+2)}{2}\right)-(i+j)}, \dots,x_1^{k}x_3^{d-k},\\
x_1^{(k-1)} x_2^{d-(k-1)}, x_1^{(k-1)} x_2^{(d-1)-(k-1)}x_3,\dots, x_1^{(k-1)} x_3^{d-(k-1)}, \\
\hskip 2.5 true cm \vdots \\
x_1^{\ell+1} x_2^{(d-1)-\ell}, x_1^{\ell+1} x_2^{(d-2)-\ell}x_3,
\dots, x_1^{\ell+1} x_3^{(d-1)-\ell}\\
x_1^\ell x_2^{d-\ell}, 
\end{array}
$$
and thus
$$
\ds \beta_{2,d+2}=\sum_{T\in \GG(I)_d} \binom{m(T)-1}{2}=j+\ell.
$$

\item $i+j=(k-1)d-\frac{(k-1)k}{2}$ and $k=\ell+1,\dots,d$.

Then the first monomial of $I_d-R_1I_{d-1}$ is
$$
x_1^{(k-1)} x_2^{d-(k-1)},
$$
and hence we have $(j+k)$-generators in $I_d$ as follows:
$$
\begin{array}{llllllllllllllll}
x_1^{(k-1)} x_2^{d-(k-1)}, x_1^{(k-1)} x_2^{(d-1)-(k-1)}x_3,\dots, x_1^{(k-1)} x_3^{d-(k-1)}, \\
\hskip 2.5 true cm \vdots \\
x_1^{\ell+1} x_2^{(d-1)-\ell}, x_1^{\ell+1} x_2^{(d-2)-\ell}x_3,
\dots, x_1^{\ell+1} x_3^{(d-1)-\ell}\\
x_1^\ell x_2^{d-\ell},  
\end{array}
$$

and thus
$$
\ds \beta_{2,d+2}=\sum_{T\in \GG(I)_d} \binom{m(T)-1}{2}=j+\ell,
$$
\end{enumerate}
\end{enumerate}
as we wished. 
\end{proof}

Now we are ready to prove the main theorem in this section.

\begin{thm} \label{T:014} 
Let $\H$ and $j$ be as in Proposition~\ref{P:014}. Then for every $-(d-1)\le i
\le d+2$, $\H$ is not level. 
\end{thm} 

\begin{proof}
By Proposition~\ref{P:006}, this theorem holds for $-(d-1) \le i \le 1$. It suffices to prove this theorem for $2 \le i \le d+2$. By Proposition~\ref{P:014}, we have that
\begin{equation} \label{EQ:003}
\begin{array}{lllllllllllll}
\beta_{1,d+2}=
\begin{cases}
2, & \text{for} \quad i=2,\dots,d, \\
3, & \text{for} \quad i=d+1,d+2,
\end{cases} 
\quad \text{ and } \quad 
\beta_{2,d+2}=
\begin{cases}
j+1, & \text{for} \quad i=2,\dots,d, \\
j+2, & \text{for} \quad i=d+1,d+2. 
\end{cases}
\end{array}
\end{equation} 
Hence if $j\ge 2$, then $\H$ is not level since $\beta_{2,d+2}>\beta_{1,d+2}$. 
It is enough, therefore,  to show the case $j=1$. 

First, assume  $i=2,3,\dots,d$. Then, by equation~\eqref{EQ:003},  we have $\beta_{1,d+2}=\beta_{2,d+2}=2$. Moreover, we see that $h_{d-1}=d+i+j=d+i+1$ and $h_d=h_{d+1}=d+i$. 

Now suppose $A=R/I,\ R=k[x_1,x_2,x_3]$, is a level algebra with $h$-vector $(h_0,h_1,\dots,h_{d-1},h_d, \linebreak 
h_{d+1})$ where $h_d=h_{d+1}$ and the ideal $I$ has $h_{d+1}^{\langle d+1\rangle}=(d+i)^{\langle d+1\rangle}=(d+i+1)$-generators in degree $d+2$. Let $J=(I_{\le d+1})$. Then the Hilbert function of $R/J$ begins
$$
h_0,h_1,\dots,h_{d-1},\overset{d\text{-th}}{d+i},\overset{(d+1)\text{-st}}{d+i},
\overset{(d+2)\text{-nd}}{d+i+1},\dots.
$$
Note that $d+i-1,d+i,d+i+1$ in degrees $d,d+1,d+2$ have the maximal growth, and so, by Theorem 3.4 in \cite{GHMS}, $R/J$ has one dimensional socle element in degree $d$. Since $R/J$ and $R/I$ agree in degree $\le d+1$, $R/I$ has such a socle element. It follows that in order for $R/I$ to be level, $I$ must have at most $(d+i)$-generators in degree $d+2$. Then both copies $R(-(d+2))$ of the last free module of the minimal free resolution of $R/I$ cannot be canceled. Therefore, the Hilbert function $\H$ cannot be level. 

Second, assume $i=d+1$. By equation~\eqref{EQ:003}, we have $\beta_{1,d+2}=\beta_{2,d+2}=3$. 

Suppose the ideal $I$ has $h_{d+1}^{\langle d+1\rangle}=(2d+1)^{\langle d+1\rangle}=(2d+2)$-generators in degree $d+2$ and let $J=(I_{\le d+1})$. Then the Hilbert function of $R/J$ begins
$$
h_0,h_1,\dots,h_{d-1},\overset{d\text{-th}}{2d+1},\overset{(d+1)\text{-st}}{2d+1},
\overset{(d+2)\text{-nd}}{2d+2},\dots.
$$
Note also that $2d,2d+1,2d+2$ in degrees $d,d+1,d+2$ have the maximal growth. Therefore, by Theorem 3.4 in \cite{GHMS} again, $R/J$ has one dimensional socle element in degree $d$, so does $R/I$ by the same argument as above. Thus three copies $R(-(d+2))$ of the last free module of the minimal free resolution of $R/I$ cannot be canceled. Therefore, the Hilbert function $\H$ cannot be level.

Finally, assume $i=d+2$. By the similar argument to the case $i=d+1$, $\H$ is not level either, as we wished. 
\end{proof}

Theorem~\ref{T:014} shows that Question~\ref{Q:012} is true if $h_d\le 2d+2$. The following example shows a case $j=1$ and $i=d+2\ (h_d=2d+2)$ of this theorem.

\begin{exmp} \label{R:058} Let $I$ be the lex-segment ideal in $R=k[x_1,x_2,x_3]$ with Hilbert function    
$$
\begin{array}{lllllllllllllllllllll}
\H & : & 1 & 3 & 6 & 10 & 15 & 21 & 17 & 16 & 16 & 0 & \ra \ .
\end{array}
$$
Note that $h_7=16=2\times 7+2=2d+2$, which satisfies the condition in Theorem~\ref{T:014}, and $j=h_6-h_7=17-16=1$. Hence
any Artinian ring with Hilbert function $\H$ cannot be level. 
\end{exmp} 

We now give another example for $i=d+3\ (h_d=2d+3)$ which does not satisfy the condition in Theorem~\ref{T:014}. 

\begin{exmp} \label{EX:059} 
Let $I$ be the lex-segment ideal of $R$ with Hilbert function
$$
\begin{array}{lllllllllllllllllllll}
\H & : & 1 & 3 & 6 & 10 & 15 & 21 & 18 & 17 & 17 & 0 & \ra.
\end{array}
$$
Note that  $h_7=17=2\times 7+3=2d+3$ and $j=18-17=1$. Hence, by Proposition~\ref{P:014}, we have $\beta_{1,d+2}=4$ and $\beta_{2,d+2}=3$, that is, $\beta_{1,d+2}>\beta_{2,d+2}$. This means that we cannot say if Hilbert function $\H$ is level only  based on shifts and Betti numbers.   In other words, for $i\le d+2$ (or $h_d \le 2d+2$), we can decide if $\H$ is {\bf\em not}  level using shifts and Betti numbers. 
\end{exmp}

We now pass to Question~\ref{Q:051} and the following corollary answers to this question for $i\le d+2$ (or $h_d \le 2d+2$). 

\begin{cor} \label{C:510}
Let $\H=\{h_i\}_{i\ge 0}$ be an O-sequence with $h_1=3$. If
$$
h_{d-1}>h_d, \quad h_d\le 2d+2, \quad \text{ and } \quad h_{d+1}\ge h_d
$$
for some degree $d$, then $\H$ is not level. 
\end{cor} 

\begin{proof} Note that, by the proof of Theorem~\ref{T:014}, any graded ring with Hilbert function 
$$
\begin{matrix}
\H' & : & h_0 & h_1 & \cdots & h_{d-1} & h_d & h_d & \ra
\end{matrix}
$$
has a socle element in degree $d-1$. 

Now let $A=\bigoplus_{i\ge 0} A_i$ be a graded ring with Hilbert function $\H$.  If $A_{d+1}=\langle f_1,f_2,\dots,f_{h_{d+1}}\rangle$ and $I=(f_{h_d+1},\dots,f_{h_{d+1}})\bigoplus_{j\ge d+2} A_j$, then a graded ring $B=A/I$ has Hilbert function
$$
\begin{matrix}
h_0 & h_1 & \cdots & h_{d-1} & h_d & h_d,
\end{matrix}
$$
and hence $B$ has a socle element in degree $d-1$ by Theorem~\ref{T:014}. Since $A_i=B_i$ for every $i\le d$, $A$ also has the same socle element in degree $d-1$ as $B$, and thus $\H$ is not level as we wished.  
\end{proof} 

\begin{exmp} Consider an O-sequence 
$$
\begin{array}{lllllllllllllllllllll}
\H & : & 1 & 3 & 6 & 10 & 14 & 18 & 17 & 16 & h_8 & \cdots .
\end{array}
$$
Then, there are only $3$ possible O-seuences such that $h_8\ge h_7=16$ since $h_8\le h_7^{\langle 7\rangle}=16^{\langle 7\rangle}=18$. By Theorem~\ref{T:014}, $\H$ is not level if $h_8=h_7=16$. The other two non-unimodal O-sequences, by Corollary~\ref{C:510},
$$
\begin{array}{lllllllllllllllllllll}
1 & 3 & 6 & 10 & 14 & 18 & 17 & 16 & 17 & \cdots \quad \text{and}\\
1 & 3 & 6 & 10 & 14 & 18 & 17 & 16 & 18 & \cdots
\end{array}
$$
cannot be level. 
\end{exmp}

\section{Degree of The Socle Elements of Graded Artinian Algebras}\label{Sec:006}

In this section, we are interested in another non-level O-sequences of codimension $3$:
\begin{equation}\label{EQ:014}
\begin{matrix}
\H & : & h_0 & h_1 & \cdots & h_{d-1} & h_d & \cdots & \overset{(d+s-1)\text{-st}}{h_d}
 & h_{d+s}
\end{matrix}
\end{equation}
where $s\ge 2$ and $h_d<h_{d+s}$. In particular, we shall prove that some graded algebra with Hilbert function $\H$ of codimension $3$ in equation~\eqref{EQ:014} has a socle element in degree $d+s-2$, and hence cannot be level.  

First, we recall the definitions of type vectors and $k$-configurations in $\P^n$. 

\begin{defn}[$n${\bf -type vectors}, Definition 2.1, {\cite{GHS:1}}] \label{D:201}
\begin{itemize}
\item[1)] A {\em $0$-type vector} will be defined to be ${\mathcal T} = 1$.  It is the only {\em $0$-type vector}.  We shall define $\alpha ({\mathcal T}) = -1$ and $\sigma ({\mathcal T}) = 1$.

\item[2)] A {\em $1$-type vector} is a vector of the form ${\mathcal T} =(d)$ where $d\geq 1$ is a positive integer.  For such a vector we define $\alpha({\mathcal T}) = d = \sigma ({\mathcal T})$.

\item[3)] A {\em $2$-type vector}, ${\mathcal T}$, is
$$
{\mathcal T} = ( (d_1), (d_2), \ldots , (d_m))
$$
where $m \geq 1$, the $(d_i)$ are {\em $1$-type vectors}.  We also insist that $\sigma (d_i)=d_i < \alpha (d_{i+1})=d_{i+1}$.
\end{itemize}

For such a $\mathcal T$ we define $\alpha({\mathcal T}) = m$ and $\sigma({\mathcal T}) =  \sigma( (d_m) ) = d_m$.  
 
Clearly, $\alpha({\mathcal T}) \leq \sigma({\mathcal T})$ with equality if and only if ${\mathcal T} = ((1),(2),\ldots , (m))$. For simplicity in the notation we usually rewrite the {\em $2$-type vector} $ ((d_1), \ldots, (d_m) ) \hbox{ as }  (d_1, \ldots , d_m) \ . $

\begin{itemize}
\item[4)] Now let $n \geq 3$.  An {\em $n$-type vector}, $\mathcal T$, is an ordered collection of {\em $(n-1)$-type vectors}, ${\mathcal T}_1, \ldots , {\mathcal T}_s$, i.e.
$$
{\mathcal T} = ({\mathcal T}_1, \ldots , {\mathcal T}_s)
$$
for which $\sigma({\mathcal T}_i) < \alpha({\mathcal T}_{i+1})$ for $i = 1, \ldots , s-1$. 
\end{itemize}
\ni
For such a $\mathcal T$ we define $\alpha({\mathcal T}) = s$ and $\sigma({\mathcal T}) = \sigma ({\mathcal T}_s)$.
\end{defn}

\def\proj{\mathbb P}
\def\cal{\mathcal}

\begin{defn}[{{\bf $k$-configuration in $\proj^n$}}, Definition 4.1, {\cite{GHS:1}}] 
\hphantom{aaaaaaaaaaaaaaaaaaaaaaaaaaaaaa}\hskip 5true cm
\begin{itemize}
\item[${\cal S}_0$:] The only element in ${\cal S}_0$ is $\H:= 1 \ \rightarrow$.  It is the Hilbert function of $\proj ^0$, which is a single point. That is the only $k$-configuration in $\proj ^0$.

\item[${\cal S}_1$:] Let $\H \in {\cal S}_1$.  Then $\chi_1(\H) = {\cal T} = (e)$ where $e \geq 1$.  We associate to $\H$ any set of $e$ distinct points in $\proj ^1$.  Clearly any set of $e$ distinct points in $\proj ^1$ has Hilbert function $\H$.

A set of $e$ distinct points in $\proj ^1$ will be called a {\em k-configuration in $\proj ^1$ of type ${\cal T} = (e)$}.

\item[${\cal S}_2$:]
Let $\H \in {\cal S}_2$ and let ${\cal T} = ( (e_1), \ldots , (e_r) ) = \chi_2(\H)$, where ${\cal T}_i = (e_i)$ is a {\em $1$-type vector}.  Choose $r$ distinct $\proj ^1$'s in $\proj ^2$ i.e. lines in $\proj ^2$, and label them ${\mathbb L}_1, \ldots, {\mathbb L}_r$.  By induction we choose, on ${\mathbb L}_i$, a {\em k-configuration} in $\proj ^1$, call it ${\mathbb X}_i$, of type ${\cal T}_i = (e_i)$ -- each {\em k-configuration} chosen so that no point of ${\mathbb L}_i$ contains any point of ${\mathbb X}_j$ for $j<i$.

The set ${\mathbb X} = \bigcup{\mathbb X}_i$ is called a {\em  k-configuration  in} $\proj ^2$  {\em of type} $\cal T$.

\item[${\cal S}_n$,] $(n>2)$:
Now suppose that we have defined a {\em k-configuration} of {\em Type} $\widetilde{\cal T}\ \in \proj ^{n-1}$, where $\widetilde{\cal T}$ is an {\em $(n-1)$-type vector} associated to $G \in {\cal S}_{n-1}$.

Let $\H \in {\cal S}_n$ and suppose that $\chi _n(\H) = {\cal T} = ({\cal T}_1, \ldots , {\cal T}_r)$ where the ${\cal T}_i$ are {\em 
$(n-1)$-type vectors}.  Then $\rho_{n-1}({\cal T}_i) = \H_i$ and $\H_i \in {\cal S}_{n-1}$.

Consider ${\mathbb H}_1, \ldots , {\mathbb H}_r$ distinct hyperplanes in $\proj ^n$ and let ${\mathbb X}_i$ be a {\em k-configuration} in ${\mathbb H}_i$ of type ${\cal T}_i$ such that ${\mathbb H}_i$ does not contain any point of ${\mathbb X}_j$ for any $j < i$.

The set ${\mathbb X} = \bigcup{\mathbb X}_i$ is called a {\em k-configuration in} $\proj ^n$ {\em of type} $\cal T$.
\end{itemize}
\end{defn}

Now we shall introduce some non-level O-sequences based on type vectors. 

\begin{rem}\label{R:061} 
\begin{itemize} 
\item[(a)] Let $\X$ be a $k$-configuration in $\P^2$ of type $\T=(d_1,\dots,d_\alpha)$ with $d_{i+1}-d_i\ge 3$ for some $i=1,\dots,\alpha-1$.   
Since $\X$ is a $k$-configuration in $\P^2$ of type $\T=(d_1,\dots,d_\alpha)$, we have the minimal free resolution of $R/I_\X$ is
$$
\begin{array}{llllllllll}
0 & \ra  R(-(d_1+\alpha))\oplus\cdots\oplus R(-(d_i+\alpha-i+1)) \oplus\cdots\oplus R(-(d_\alpha+1)) \\
  & \ra  R(-\alpha)\oplus R(-(d_1+\alpha-1))\oplus\cdots\oplus R(-(d_i+\alpha-i)) \oplus\cdots\oplus R(-d_\alpha) \\
  & \ra  R \ra  R/I_\X \ra \ 0
\end{array}
$$
by Theorem 2.6 in~\cite{GPS}. Since $d_{i+1}-d_i\ge 3$, we have that $d_i+\alpha-i+1<d_{i+1}+\alpha-(i+1)$, which means that $R(-(d_i+\alpha-i+1))$ of the last free module cannot be canceled. Hence the Hilbert function $\H_\X$ is not level.

\item[(b)] Let $\X$ be a $k$-configuration in $\P^3$ of type $\mathcal T=(\mathcal T_1,\dots,\mathcal T_\alpha)$ and let $\FF_\X$ be the minimal free resolution of the coordinate ring of $\X$. If either $\mathcal T_i$ is the $2$-type vector as in this remark (a) or $\sigma(\mathcal T_i)+2<\alpha(\mathcal T_{i+1})$ for some $i=1,\dots,\alpha-1$, then the Hilbert function $\H_\X$ is not level. To show this, we shall use the same notation as in  Theorem 3.2, \cite{GS:1}.  If  $\T_i$ is a $2$-type vector as in (a) for some $i$, then $\H_\X$ is obviously not level by the same idea as in (a). Now assume that $\sigma(\mathcal T_i)+2<\alpha(\mathcal T_{i+1})$. Then we have
$$
\begin{array}{rllllllllllll}
\varepsilon_i+2+\bar d_{i\alpha(\mathcal T_i)}
& = & \alpha(\mathcal T_i)-i+\alpha+2+\sigma(\T_i)-\alpha(\T_i) \\
& = & \sigma(\T_i)-i+\alpha+2, \quad \text{and} \\
\varepsilon_{i+1}+1
& = & \alpha(\T_{i+1})-(i+1)+\alpha+1\\
& = & \alpha(\T_{i+1})-i+\alpha.
\end{array}
$$
Since $\sigma(\T_i)+2<\alpha(\T_{i+1})$, we see that 
$\varepsilon_i+2+\bar d_{i\alpha(\mathcal T_i)}
< \varepsilon_{i+1}+1$, and hence $R(-(\varepsilon_i+2+\bar d_{i\alpha(\mathcal T_i)}))$ in the last free module of $\FF_\X$ cannot be calceled, that is, the Hilbert function $\H_\X$ is not level.  
 
\item[(c)] Let $\T=(\T_1,\dots,\T_\alpha)$ be the $3$-type vector with $\T_i=(d_{i1},\dots,d_{i\alpha(\T_i)})$ for every $i$. If $\alpha(\T_i)=\sigma(\T_{i-1})+1$ and $d_{i1}\ge 3$ for some $i$, then the Hilbert function of a $k$-configuration $\X$ in $\P^3$ of type $\T$ is not level. To show this, we shall use the same notation as in  Theorem 3.2, \cite{GS:1} again. Then we have
$$
\begin{array}{rllllllllllll}
\varepsilon_{i-1}+2+\bar d_{i-1\alpha(\T_{i-1})}
& = & \alpha(\mathcal T_{i-1})-(i-1)+\alpha+2+d_{i-1\alpha(\T_{i-1})}-\alpha(\T_{i-1}) \\
& = & d_{i-1\alpha(\T_{i-1})}+\alpha-i+3, \quad \text{and} \\
\varepsilon_{i}+1
& = & \alpha(\T_{i})+\alpha-i+1.
\end{array}
$$
Since $\alpha(\T_{i})=\sigma(\T_{i-1})+1=d_{i-1\alpha(\T_{i-1})}+1$, we see that 
$\varepsilon_{i}+1<\varepsilon_{i-1}+2+\bar d_{i-1\alpha(\T_{i-1})}$. Moreover, we have 
$$
\begin{array}{llllllllllll}
\varepsilon_i+1+\bar d_{i1}
& = & \alpha(\T_{i})+\alpha-i+1+d_{i1}-1 \\
& = & \alpha(\T_{i})+\alpha-i+d_{i1} \\
& > & d_{i-1\alpha(\T_{i-1})}+\alpha-i+3 \\
& = & \varepsilon_{i-1}+2+\bar d_{i-1\alpha(\T_{i-1})}.
\end{array}
$$
In other words, 
$$
\varepsilon_{i}+1<\varepsilon_{i-1}+2+\bar d_{i-1\alpha(\T_{i-1})}< \varepsilon_i+1+\bar d_{i1},
$$
and hence $R(-(\varepsilon_{i-1}+2+\bar d_{i-1\alpha(\T_{i-1})}))$ in the last free module of $\FF_\X$ cannot be canceled, that is, the Hilbert function $\H_\X$ is not level.
\end{itemize}
\end{rem}

Now we are ready to discuss about the degree of the socle elements of some Artinian algebra with Hilbert function $\H$ in equation~\eqref{EQ:014}. 

\begin{thm} \label{T:062} Let $\H$ be as in equation~\eqref{EQ:014} and $\T=(\T_1,\dots,\T_\alpha)$ be the $3$-type vector corresponding to the Hilbert function whose first difference is $\H$. 
If $h_d=\cdots=h_{d+s-1}=d+s+(i-1)$ and $h_{d+s}=d+s+i$ where $1\le i\le \alpha(\T_{\alpha-1})$, then 
$$
\begin{array}{rlllllllllllllll}
\T_\alpha
& = & (1,2,\dots,d+s,d+s+1), \\[2ex] 
\T_{\alpha-1}
& = & 
\begin{cases}
(\dots,d+s-2),  & \text{for } i=1,\\
(\dots,d+s-(i+1),d+s-(i-2),\dots,d+s),  
 & \text{for } i=2,\dots,\alpha(\T_{\alpha-1}).
\end{cases}
\end{array}
$$ 
In particular, the O-sequence $\H$ is not level and any Artinian graded algebra with Hilbert function $\H$ has a socle element in degree $d+s-2$.
\end{thm}

\begin{proof} It suffices to prove this theorem for $i=1$ and $2$, respectively, since we can use the same argument for the rest of the cases $i\ge 3$ as for $i=2$.  

\begin{itemize}
\item[] {\em Case 1.} If $i=1$, that is, $h_d=\cdots=h_{d+s-1}=d+s$ and $h_{d+s}=d+s+1$, then from the following equation, 

\vskip 2pc

\phantom{a}

\begin{figure}[ht]
\vskip -3pc
\centering
\begin{picture}(220,-200)
\put(-20,-60){$\begin{matrix} 
\H & : & h_0 & h_1 & h_2 & \cdots & \overset{(d+s-2)\text{-nd}}{h_d} & \overset{(d+s-1)\text{-st}}{h_d}     
              & h_{d+s} \\
          &   & 1 & 2 & 3 & \cdots  & d+s-1   & d+s & d+s+1 \\[2ex]
          &   &   & 1 & a & \cdots  & 1 & 0 & 0 \\
          &   &   & 1 &   & \cdots & 1 &   \\[2ex] 
          &   &   &   &   & \cdots      
\end{matrix}$}
\put(-25,-45){\line(1,0){290}}
\put(-25,-92){\line(1,0){290}}
\end{picture}
\vskip 3.5true cm
\end{figure}
\ni
where $1\le a \le 3$, we have that
$$
\begin{array}{rlllllllllllllll}
\T_\alpha
& = & (1,2,3,\dots,d+s,d+s+1) \\
\T_{\alpha-1}
& = & (\dots,d+s-2). 
\end{array}
$$
In what follows
$$
\alpha(\T_\alpha)=d+s+1=(d+s-2)+3=\sigma(\T_{\alpha-2})+3 > \sigma(\T_{\alpha-1})+2,
$$
and hence, by Remark~\ref{R:061}~(b), $\H$ is not level. 

Recall that $\FF_\X$ is the minimal free resolution of a coordinate ring $R/I_\X$ of a $k$-configuration $\X$ in $\P^3$ with Hilbert function $\mathbf G$ such that  $\Delta\mathbf G=\H$ and we shall use the same notation as in Theorem 3.2 in \cite{GS:1} for the rest of the proof. Since  the non-cancelable shift of the last free module of $\FF_\X$ is
$$
\begin{array}{lllllllllllllll}
  & \varepsilon_{\alpha-1}+2+\bar d_{\alpha-1\alpha(\T_{\alpha-1})}\\
= & [\alpha(\T_{\alpha-1})-(\alpha-1)+\alpha]+2+d_{\alpha-1\alpha(\T_{\alpha-1})}
-\alpha(\T_{\alpha-1}) \\
= & d_{\alpha-1\alpha(\T_{\alpha-1})}+3 \\
= & (d+s-2)+3,
\end{array}
$$ 
any algebra with Hilbert function $\H$ has a socle element in degree $d+s-2$. 

\item[] {\em Case 2.} If $i=2$, that is, $h_d=\cdots =h_{d+s-1}=d+s+1$ and $h_{d+s}=d+s+2$, then from the following equation  

\begin{figure}[ht]
\centering
\begin{picture}(220,0)
\put(-20,-60){$\begin{matrix} 
\H & : & h_0 & h_1 & h_2 & \cdots & h_d & \overset{(d+s-1)\text{-th}}{h_d}     
              & h_{d+s} \\
          &   & 1 & 2 & 3 & \cdots  & d+s-1   & d+s & d+s+1 \\[2ex]
          &   &   & 1 & a & \cdots  & 2 & 1 & 1 \\
          &   &   & 1 & 2 & \cdots & 2 &  1 & 1 \\[2ex] 
          &   &   &   &   & \cdots      
\end{matrix}$}
\put(-25,-45){\line(1,0){290}}
\put(-25,-92){\line(1,0){290}}
\end{picture}
\vskip 3.5true cm
\end{figure}
\ni
where $2\le a \le 3$, we have that
$$
\begin{array}{rllllllllllllll}
\T_\alpha
& = & (1,2,3,\dots, d+s,d+s+1) \\
\T_{\alpha-1}
& = & (\dots,d+s-3,d+s).
\end{array}
$$
Since the difference of the last two 2-type vectors of $\T_{\alpha-1}$ is $3$, by Remark~\ref{R:061}~(a), any Artinian algebra with Hilbert function $\H$ is not level. Furthermore, the non-cancelable shift of the last free module of $\FF_\X$ is
$$
\begin{array}{lllllllllllllll}
  & \varepsilon_{\alpha-1}+2+\bar d_{\alpha-1(\alpha(\T_{\alpha-1})-1)}\\
= & [\alpha(\T_{\alpha-1})-(\alpha-1)+\alpha]+2+d_{\alpha-1(\alpha(\T_{\alpha-1})-1)}
-(\alpha(\T_{\alpha-1})-1) \\
= & d_{\alpha-1(\alpha(\T_{\alpha-1})-1)}+4 \\
= & (d+s-2)+3, 
\end{array}
$$
and thus any Artinian algebra with Hilbert function $\H$ has a socle element in degree $d+s-2$. 
\end{itemize}
By continuing this process for $3 \le i\le \alpha(\T_{\alpha-1})$, we complete the proof, as we wished. 
\end{proof} 

\begin{rem}  Using the same argument as in the proof of Corollary~\ref{C:510}, Theorem~\ref{T:062} holds when $h_{d+s}>h_d$, in general. 
\end{rem}

\def\G{\mathbf G}

\begin{exmp} \label{EX:063}
Consider an O-sequence $\begin{matrix} 
\H & : & 1 & 3 & 6 & 8 & 9 & 9 & 9 & 10 \end{matrix}.$ Then $d=4$ and $s=3$, and so $h_{d+s}=10=7+3$, that is, $i=3$. Note that $\alpha(\T_2)=4$. Applying Theorem~\ref{T:062} to this case,  we conclude any graded Artinian algebra with Hilbert function $\H$ is not level and has a socle element in degree $d+s-2=5$.
\end{exmp}

\end{document}